\documentclass[11pt,oneside,a4paper]{article}
\usepackage{syntonly}
\usepackage{geometry}
\usepackage{mathrsfs}
\usepackage{makeidx}
\makeindex
\usepackage{amsmath}
\usepackage{amssymb}
\usepackage{amsmath}
\usepackage{amsfonts}
\usepackage{amssymb}
\usepackage{esvect}
\usepackage{algorithm}
\usepackage{bm}
\usepackage{tikz}
\usepackage{cleveref}

\usepackage{graphicx}
\usepackage{epstopdf}

\newtheorem{theorem}{Theorem} 
\newtheorem{lemma}{Lemma}
\newtheorem{proposition}{Proposition}
\newtheorem{corollary}{Corollary}

\newcommand{\qed}{\nobreak \ifvmode \relax \else
      \ifdim\lastskip<1.5em \hskip-\lastskip
      \hskip1.5em plus0em minus0.5em \fi \nobreak
      \vrule height0.30em width0.4em depth0.25em\fi}

 \author{Safari Mukeru\\
\footnotesize{\em Department of Decision Sciences}\\ 
\footnotesize{University of South Africa, P. O. Box 392, Pretoria, 0003. South Africa}\\
\footnotesize{e-mail: mukers@unisa.ac.za}}

\title{{On the convergence of series of dependent  random variables}}

\date{}

\begin{document}

\maketitle

\pagenumbering{arabic}

\begin{abstract}
Given a sequence $(X_n)$ of symmetrical random variables taking values in a Hilbert space, an interesting open problem is to determine the conditions under which  the series $\sum_{n=1}^\infty X_n$ is almost surely convergent. For independent random variables, it is well-known that if $\sum_{n=1}^\infty \mathbb{E}(\|X_n\|^2) <\infty$, then $\sum_{n=1}^\infty X_n$ converges almost surely. This has been extended to some cases of dependent variables (namely negatively associated random variables) but in the general setting of dependent variables, the problem remains open. This paper considers the case where each variable $X_n$ is given as a linear combination $a_{n,1}Z_1+ \ldots +a_{n,n}Z_n$ where  $(Z_n)$ is a sequence of independent symmetrical random variables of unit variance and $(a_{n,k})$ are constants.  For Gaussian random variables, this is the general setting. We obtain a sufficient condition for the almost sure convergence of  $\sum_{n=1}^\infty X_n$ which is also sufficient for the almost sure convergence of  $\sum_{n=1}^\infty \pm X_n$ for all (non-random) changes of sign. 
The result is based on an important bound of the mean of the random variable 
$\sup(\|X_1 + \ldots +X_k\|: 1\leq k \leq n)$ 
which extends the classical L\'evy's inequality and has some independent interest. 
\end{abstract}
{\bf Key words:} random series, almost sure convergence, $L^2$--convergence, Hilbert spaces. \\ 
60B12, 60G50, 40A05.

\section{Introduction}
Given a Hilbert space $\mathbb{H}$, a probability space $(\Omega, \mathcal{F}, \mathbb{P})$ and a sequence of random variables $(X_n)$ taking values in $\mathbb{H}$, under which conditions is the random series $S\equiv \sum_{n=1}^\infty X_n $ almost surely convergent (in the sense that the sequence $(S_n)$ of partial sums $S_n = X_1+ \ldots+X_n$ converges in norm almost surely)? 
In particular, for a sequence of vectors $(u_n)$ in $\mathbb{H}$ and a sequence of real valued random variables $X_n: \Omega \to \mathbb{R}$, $n =1,2, \ldots,$ under which conditions on the vectors $(u_n)$ is the random series 
 $S(\omega) \equiv \sum_{n=1}^\infty X_n(\omega) u_n$ 
almost surely convergent?  This problem has a long and rich history and  is associated with the names of Kolmogorov, Rademacher, Zygmund, Paley, L\'evy, It\^o, Kahane, Pisier, etc. It is also related to the problem of re-arrangements of series (see for example, Kvaratskhelia \cite{Kvaratskhelia} and  Levental, Mandrekar and Chobonyan \cite{Levental}) and various limit theorems such as the central limit theorem, the law of large numbers, etc. 
A classical result (due originally to Rademacher) says that if $(\xi_n)$ is the sequence of independent and identically distributed random variables with $\mathbb{P}\{\xi_n = 1\} = \mathbb{P}\{\xi_n = -1\} = \frac{1}{2}$ and $(u_n)$ is a sequence of vectors in a Hilbert space $\mathbb{H}$ such that $\sum_{n=1}^\infty \|u_n\|^2 < \infty$, then the series $\sum_{n=1}^\infty \xi_n u_n$ converges almost surely. This has been generalised to sequences of symmetrical independent random variables (see for example Kahane \cite[Theorem 2, p 30]{Kahane_1985}): for any sequence $(X_n)$ of symmetrical independent random variables 
if $\sum_{n=1}^\infty \mathbb{E}(\|X_n\|^2) < \infty$, then the series $\sum_{n=1}^\infty X_n$ converges almost surely. 

Some authors have considered the problem of convergence of series of {\it dependent} random variables in the specific case of sequences of  {\it associated} and {\it negatively associated} random variables introduced by Esary, Proschan and Walkup \cite{Esary} and Joag-Dev and Proschan \cite{Joag-Dev}.   
{Matu\l{}a} \cite{Matula} proved that if $(X_n)$ is a sequence of negatively associated real-valued random variables with finite second moments and zero means and if $\sum_{n=1}^\infty \mathbb{E}(\|X_n\|^2) < \infty$, then the series $\sum_{n=1}^\infty X_n$ converges almost surely. This result was extended to negatively associated random variables taking values in Hilbert spaces by Ko, Kim and Han \cite{Ko} (see also Wu and Jiang \cite{Wu_Jiang} for further properties of negatively associated random variables). However, in the general setting, the almost sure convergence of random series of dependent random variables remains so far an open problem. 

The study of series of dependent random variables is a difficult problem because most classical tools that are useful in the case of independent variables are no longer available in the general setting. Important examples are the classical L\'evy's inequality and the Kolmogorov maximum inequality  
which play a key role in analysing random series of independent random variables but they do not hold for the general case of dependent random variables. 

In this paper, we consider the general setting of random series $\sum_{n=1}^\infty X_n$ such that each random variable $X_n$ is a linear combination $a_{n,1}Z_1+ \ldots +a_{n,n}Z_n$ where  $(Z_n)$ is a fixed sequence of independent and symmetrical random variables of unit variance taking values in a Hilbert space $\mathbb{H}$ and $(a_{n,k})$ are complex numbers. 
It is a classical result that from any sequence $(X_k)$ of symmetrical  random variables (of finite second moment) in a Hilbert space, one can construct a sequence of symmetrical random variables $(Z_k)$ taking values in $\mathbb{H}$ that is orthonormal (in the sense that $\mathbb{E}(\langle Z_k,  Z_j\rangle) = 0$ for  $k\ne j$ and $\mathbb{E}(\|Z_k\|^2) = 1$) and such that each $X_n$ is a linear combination of $Z_1, Z_2, \ldots, Z_n$ (see for example, Shiryaev \cite[pp 318--323]{Shiryaev}). In the particular case where $(X_n)$ are all Gaussian random variables,  then $(Z_n)$ are also Gaussian and hence independent. The situation considered here is therefore very broad. In this paper we obtain a simple and  explicit analytical condition which guarantees the almost sure convergence of the series $\sum_{n=1}^\infty X_n$. 
This is also extended to the case where each coefficient $a_{n,k}$ is a random variable taking value in the complex plane under the additional requirement that $a_{n,k}$ is measurable with respect to the $\sigma$-algebra spanned  by $Z_1, Z_2, \ldots, Z_{k-1}$.

\section{Main results}
A naive solution to the problem of almost sure convergence of the series $\sum_{n=1}^\infty X_n$ where $X_n = \sum_{k=1}^n a_{n,k} Z_k$  would be to write 
$$\sum_{n=1}^\infty X_n = \sum_{n=1}^\infty \left(\sum_{k=1}^n a_{n,k} Z_k\right)$$
and then {\it re-arrange} the terms of the series in order to obtain 
     \begin{eqnarray} \label{s423de234sdw}
 \sum_{n=1}^\infty X_n = \sum_{k=1}^\infty \left(\sum_{n=1}^\infty a_{n, k}\right) Z_k.
\end{eqnarray} 
Then one would obtain that 
         $$\sum_{k=1}^\infty \left|\sum_{n=1}^\infty a_{n, k}\right|^2 < \infty$$  guarantees the almost sure convergence of the initial series $\sum_{n=1}^\infty X_n$. The problem with this argument is that the general setting considered in this paper does not justify the {\it re-arrangement} of the terms  in order to obtain (\ref{s423de234sdw}).  We shall however obtain a sufficient condition on the coefficients $a_{n,k}$ for the convergence of the series $\sum_{n=1}^\infty X_n$. \\
Our main results are as following. 
\begin{theorem}
\label{th01} 
Let $(Z_n)$ be a sequence of independent symmetrical random variables with unit variance taking values in a Hilbert space $\mathbb{H}$ and let $(X_n)$ be a sequence given by \\ $X_n = a_{n,1} Z_1~+ \ldots a_{n,n}Z_n$  where $(a_{n,k})$, $1\leq k \leq n < \infty$, is a sequence of complex numbers.
If  \begin{eqnarray} \label{sder434e}
\sum_{n=1}^\infty \left(\sum_{k=1}^\infty |a_{n+k-1,k}|^2\right)^{1/2} < \infty,
\end{eqnarray}
 then
the series $\sum_{n=1}^\infty X_n$ converges in $\mathbb{H}$ almost surely. Moreover for each (non-random) change of sign, the series $\sum_{n=1}^\infty \pm X_n$ also converges in $\mathbb{H}$ almost surely. In particular if $(\epsilon_n)$ is a change of sign independent of the sequence $(Z_n)$, then the series $\sum_{n=1}^\infty \epsilon_n  X_n$ converges almost surely.  
\end{theorem}  
In the case where the random variables $(X_n)$ are independent (which implies that $a_{n,k} = 0$ for $n\ne k$), then (\ref{sder434e}) reduces to the classical sufficient condition $\sum_{n=1}^\infty \mathbb{E}(\|X_n\|^2) <\infty$ for the almost sure convergence of $\sum_{n=1}^\infty X_n$. This condition is also necessary under the additional condition that $\mathbb{E}(\|X_n\|^4) < C\ \mathbb{E}(\|X_n\|^2) < \infty$ for some fixed constant $C>0$ (Kahane \cite[p 31]{Kahane_1985}, Theorem 4).  
 Theorem \ref{th01} is a natural extension of the classical result on the convergence of random series of independent random variables.  The other extreme case is where all the $X_n$ are collinear random variables in $\mathbb{H}$ which implies  $a_{n,k} = 0$ for all $k\geq 2$. Then condition (\ref{sder434e}) reduces to $\sum_{n=1}^\infty |a_{n,1}| < \infty$ and hence $\sum_{n=1}^\infty \mathbb{E}(|X_n|) < \infty$. This is also obviously a necessary condition for the almost sure convergence of all the changes of sign series $\sum_{n=1}^\infty \pm X_n$ (in the case of collinearity).  Since our sufficient condition (\ref{sder434e}) is also a necessary condition in these two extreme cases, it would be interesting  to know if that can be extended to the general setting. 
 
Another particular case is where there is $m$ fixed such that each $X_n$ is a linear combination of $Z_{n}, Z_{n-1}, \ldots, Z_{n-m}$ only. (This is equivalent to say that $a_{n,k} = 0$ for all $k \leq n-m$). Here the condition $\sum_{n=1}^\infty \mathbb{E}(\|X_n\|^2) < \infty$ is sufficient for the almost sure convergence of $\sum_{n=1}^\infty \pm X_n$. An example is obtained from the  fractional Gaussian noise of index 0: Given a sequence $(Z_n)$  of real-valued independent standard Gaussian  random variables, set for $n=1,2,\ldots$,
$$\Delta_n =-\left(\frac{n-1}{2n}\right)^{1/2} Z_{n-1} + \left(\frac{n+1}{2n}\right)^{1/2} Z_n,\,\,\,\,\, Z_0 = 0.$$
(The sequence $(\Delta_n)$ is a model of the fractional Gaussian noise with index 0.) 
Clearly  $\mathbb{E}(\Delta_n^2) = 1$, $\mathbb{E}(\Delta_n \Delta_k) = -1/2$ for $|n-k| = 1$ and $\mathbb{E}(\Delta_n \Delta_k) = 0$ for $|n-k|> 1.$     
For any sequence $(u_n)$ of vectors in a Hilbert space $\mathbb{H}$, if $\sum_{n=1}^\infty \|u_n\|^2 < \infty$, then the series $\sum_{n=1}^\infty \pm \Delta_n u_n$ converges almost surely in $\mathbb{H}$ (for all non-random changes of sign). 
 
 Note that the sum in (\ref{sder434e}) is obtained as follows: take the infinite triangular matrix $(a_{n,k})$, $1\leq k \leq n < \infty$ and consider its diagonals $d_n$, ($n\geq 1$) given by 
\begin{eqnarray*}
d_1 & = & (a_{1,1}, a_{2,2}, a_{3,3},\ldots) = (a_{k,k})_{k\geq 1}\\ 
d_2 & = & (a_{2,1}, a_{3,2}, a_{4,3}, \ldots) = (a_{k+1, k})_{k\geq 1} \\
\vdots\\
d_n & = & (a_{n, 1}, a_{n+1, 2}, a_{n+2, 3},\ldots) = (a_{n+k-1, k})_{k\geq 1}
\end{eqnarray*}
Then the sum in (\ref{sder434e}) is given by
    $$\sum_{n=1}^\infty \left(\sum_{k=1}^\infty |a_{n+k-1,k}|^2\right)^{1/2} = \sum_{n=1}^\infty \|d_n\|_{\ell^2}.$$
In the literature, more often, random series in Hilbert spaces are of the form  $\sum_{n=1}^\infty X_n u_n$ where $(X_n)$ is a sequence of real or complex random variables and $(u_n)$ is a sequence of vectors in a Hilbert spaces. 
In that case,  Theorem \ref{th01} yields immediately the following result.
\begin{corollary}
\label{th0New} 
Let $(Z_n)$ be a sequence of real or complex independent symmetrical random variables with unit variance,  $(X_n)$ be a sequence given by \\ $X_n = a_{n,1} Z_1~+ \ldots a_{n,n}Z_n$  where $(a_{n,k})$, $1\leq k \leq n < \infty$, is a sequence of complex numbers and $(u_n)$ be a sequence of vectors in a Hilbert space $\mathbb{H}$. 
If  \begin{eqnarray} \label{sdsde34er434e}
\sum_{n=1}^\infty \left(\sum_{k=1}^\infty |a_{n+k-1,k}|^2 \|u_{n+k-1}\|^2\right)^{1/2} < \infty,
\end{eqnarray}
 then
the series $\sum_{n=1}^\infty X_n u_n$ converges in $\mathbb{H}$ almost surely. Moreover for each (non-random) change of sign, the series $\sum_{n=1}^\infty \pm X_n u_n$ also converges in $\mathbb{H}$ almost surely. 
\end{corollary}

The proof of Theorem~\ref{th01} is based on the following important inequality which is an extension of the classical L\'evy inequality.
  \begin{lemma} \label{lemma11}
 Let   $(Z_n)$ be a sequence of independent and symmetrical random variables with unit variance in a  Hilbert space $\mathbb{H}$ and  $(X_n)$ a sequence defined from a fixed scalar sequence $(a_{n,k})$ by  $X_n = a_{n,1} Z_1+ \ldots + a_{n,n}Z_n$. Then for all $N \geq 1$,
   \begin{eqnarray} \label{smewqweq234}
\mathbb{E}\left(\sup_{1\leq n \leq N}\|X_1 + X_2 + \ldots+ X_n\|\right) \leq 2
\sum_{n=1}^N \left(\sum_{k=1}^{N-n+1} |a_{n+k-1,k}|^2\right)^{1/2}.
\end{eqnarray}
  \end{lemma}  

One may ask what can be said in the case where the coefficients $(a_{n,k})$ are allowed to be also random variables rather than being restricted to constant numbers. The general case remains unsolved but the case where each $a_{n,k}$ is allowed to be measurable with respect to the $\sigma$-algebra $\sigma(Z_1, Z_2, \ldots, Z_{k-1})$ for each $n$ is given in the following result which is an interesting generalisation of Theorem \ref{th01}.
\begin{theorem} \label{thmart}
Let $(Z_n)$ be a sequence of independent symmetrical random variables with unit variance taking values in a Hilbert space $\mathbb{H}$ and let $(X_n)$ be the sequence given by \\ $X_n = a_{n,1} Z_1~+ \ldots a_{n,n}Z_n$  where $a_{n,k}: \Omega \to \mathbb{C}$ are random variables such that for each  $1\leq k \leq n < \infty$,  $a_{n,k}$ is measurable with respect to the $\sigma$-algebra $\sigma(Z_1, Z_2, \ldots, Z_{k-1})$ spanned by $Z_1, Z_2, \ldots, Z_{k-1}$. If  \begin{eqnarray} \label{sqw21der434e}
\sum_{n=1}^\infty \left(\sum_{k=1}^\infty \mathbb{E}(|a_{n+k-1,k}|^2)\right)^{1/2} < \infty
\end{eqnarray}
 then the series $\sum_{n=1}^\infty  X_n$ converges in $\mathbb{H}$ almost surely. The same holds for each change of sign $\sum_{n=1}^\infty \pm X_n$. 
\end{theorem}
The proof of Theorem \ref{thmart} is based on the following result: 
 \begin{lemma} \label{lemmamart}
 Let  $(Z_n)$ be a sequence of independent and symmetrical random variables with unit variance in a  Hilbert space $\mathbb{H}$ and let  $(X_n)$ be a sequence  defined  by  $X_n = a_{n,1} Z_1+ \ldots a_{n,n}Z_n$, where $a_{n,k}: \Omega \to \mathbb{C}$ are random variables such that for each  $a_{n,k}$ is measurable with respect to the $\sigma$-algebra $\sigma(Z_1, Z_2, \ldots, Z_{k-1})$ spanned by $Z_1, Z_2, \ldots, Z_{k-1}$. Then for each $N \geq 1$, the following inequality holds: 
   \begin{eqnarray} \label{smqw231wer}
\mathbb{E}\left(\sup_{1\leq n \leq N}\|X_1 + X_2 + \ldots+ X_n\|\right) \leq 2
\sum_{n=1}^N \left(\sum_{k=1}^{N-n+1} \mathbb{E}\left(\|a_{n+k-1,k}\|^2\right)\right)^{1/2}.
\end{eqnarray}
  \end{lemma}
  
As already discussed, in the particular case of independent random variables,  condition (\ref{sder434e}) of Theorem \ref{th01} is equivalent to $\sum_{n=1}^\infty \mathbb{E}(\|X_n\|^2) < \infty$ or equivalent the series 
$\sum_{n=1}^\infty X_n$ converges in the $L^2$--sense. Thus for symmetrical independent random variables, $L^2$--convergence implies almost sure convergence. 
It is shown in Kahane \cite[p 31]{Kahane_1985} (Theorem 4), that the converse is also true under the additional condition that  $\mathbb{E}(\|X_n\|^4) < C\ \mathbb{E}(\|X_n\|^2) < \infty$ for some fixed constant $C>0$. 
 The following result is a generalisation of that converse for dependent random variables. 

\begin{theorem} \label{eded234231s}
For any sequence $(X_n)$ of symmetrical random variables of finite second moments defined in a Hilbert space $\mathbb{H}$,  if the series $\sum_{n=1}^\infty  X_n$ converges almost surely in $\mathbb{H}$ and  there exists a constant $K>0$ such that  for all integers $m, j$,
         \begin{eqnarray} \label{cond2}
 \mathbb{E}\left(\|X_{m+1} + \ldots +  X_{m+j}\|^4\right) \leq K\ \left(\mathbb{E}\left(\|X_{m+1} + \ldots +  X_{m+j}\|^2\right)\right)^2
 \end{eqnarray}
then the series $\sum_{n=1}^\infty X_n$ also converges in $\mathbb{H}$ in the $L^2$--sense. 
\end{theorem}
This condition is clearly satisfied for symmetrical Gaussian random variables with $K = 3$.

The interplay between $L^2$--convergence and almost sure convergence of series for more general random variables can be made more precise by using the key concept of {\it stopping time} which is prominent in probability theory. We shall obtain that for any sequence $(X_n)$ of symmetrical random variables in a Hilbert space, the almost sure convergence of the series $\sum_{n=1}^\infty X_n$ can be derived from the $L^2$-convergence of another random series. 
Let $\mathscr{F}_n \equiv \sigma(X_1, X_2, \ldots, X_n)$ be the $\sigma$-algebra spanned by $X_1, X_2, \ldots, X_n$ and $\mathscr{F}_0 = \{\Omega, \emptyset\}$.  A random variable $\tau: \Omega \to \mathbb{N}$ is a stopping time if for each $k\in \mathbb{N}$, the event $\{\omega\in \Omega:\tau(\omega) = k\} \in \mathscr{F}_k$. For every stopping time $\tau$, we shall associate a sequence of random variables $(\zeta_n)$ defined by $\zeta_n(\omega) = 1$ if $\tau(\omega) \geq n$ and $\zeta_n(\omega) = -1$ if $\tau(\omega) \leq  n-1$ and consider the series $\sum_{n=1}^\infty \zeta_n X_n$. Clearly, the random variable $\zeta_n$ is measurable with respect to $\mathscr{F}_{n-1}$ since $\{\zeta_n = -1\} = \{\tau \leq n-1\} \in \mathscr{F}_{n-1}$  and $\{\zeta_n = 1\} = \{\zeta_n = -1\}^c \in \mathscr{F}_{n-1}$. \\
We have the following result:
\begin{theorem}
\label{stopth}
If for each stopping time $\tau$, the associated series $\sum_{n=1}^\infty \zeta_n X_n$ converges in the $L^2$--sense, then the series $\sum_{n=1}^\infty X_n$ converges almost surely. 
\end{theorem}         
 In the particular case where the random variables $(X_n)$ are independent, it can be readily be seen that 
   $$\mathbb{E}(\|\zeta_1 X_1 + \zeta_2X_2 + \ldots+ \zeta X_n\|^2) = \sum_{n=1}^n \mathbb{E}(\|X_n\|^2) $$
because for $k < j$,
      $$\mathbb{E}(\zeta_k \zeta_j X_k X_j) = \mathbb{E}(\zeta_k \zeta_j X_k) \mathbb{E}(X_j)  = 0$$
 since $X_j$ is independent of $X_k, \zeta_k, \zeta_{j-1}$. Then Theorem \ref{stopth} reduces to the classical result that $L^2$--convergence of random series implies almost sure convergence for symmetrical independent random variables.

\section{Proofs of the results}
\subsection{Proofs of Lemma \ref{lemma11} and Lemma \ref{lemmamart} }
We have that
\begin{eqnarray*}
X_1 & = &a_{1,1}Z_1, \\
X_2 &= &a_{2,1} Z_1 + a_{2,2}Z_2,\\
X_3& = & a_{3,1}Z_1 + a_{3,2}Z_2 + a_{3,3}Z_3,\\
\vdots \\
X_N& = & a_{N,1} Z_1 + a_{N,2}Z_2 + a_{N,3}Z_3 + \ldots+ a_{N,N}Z_N.                                                  
\end{eqnarray*}
Since $N$ is {\it finite}, we can {\it re-arrange} the terms to write
\begin{eqnarray*}
X_1 & = &a_{1,1}Z_1, \\
X_1+ X_2 &= &(a_{1,1}Z_1 + a_{2,2} Z_2) + a_{2,1}Z_1,\\
X_1+X_2+X_3& = & (a_{1,1}Z_1 + a_{2,2} Z_2 + a_{3,3}Z_3) + (a_{2,1}Z_1 + a_{3,2}Z_2) + a_{3,1}Z_1\\
\vdots \\
X_1 + X_2 +X_3+ \ldots+X_N & = & (a_{1,1}Z_1 + a_{2,2} Z_2 + a_{3,3}Z_3 + \ldots + a_{N,N}Z_N) \\
&& + (a_{2,1}Z_1 + a_{3,2}Z_2 + \ldots+ a_{N, N-1}Z_{N-1})\\
&& + (a_{3,1}Z_1 + a_{4,2}Z_2 + \ldots+ a_{N, N-2}Z_{N-2})\\
&& + \ldots \\
&&+ (a_{N-1, 1} Z_1 + a_{N, 2}Z_2) \\
&&+ a_{N,1}Z_1.                                                  
\end{eqnarray*}
Then we decompose the sequence $(X_1, X_1+X_2, \ldots, X_1+ X_2+\ldots+X_N)$ into a sum $s_{1,N} + s_{2,N} + \ldots+ s_{N,N}$ where
\begin{eqnarray*}
s_{1,N} & = &\left(a_{1,1}Z_1, a_{1,1}Z_1 + a_{2,2}Z_2, \ldots, a_{1,1}Z_1 + a_{2,2}Z_2+ \ldots + a_{N,N}Z_N\right)\\
s_{2, N} & = &\left(0, a_{2,1}Z_1 , a_{2,1}Z_1+ a_{3,2}Z_2, \ldots,  a_{2,1}Z_1+ a_{3,2}Z_2 +\ldots+ a_{N, N-1}Z_{N-1}\right)\\
s_{3, N}& = & \left(0,0, a_{3,1}Z_1, a_{3,1}Z_1+ a_{4,2}Z_2, \ldots, a_{3,1}Z_1+ a_{4,2}Z_2+ \ldots + a_{N,N-2}Z_{N-2}\right)\\
s_{N,N} &= & \left(0,0,\ldots, 0, a_{N,1}Z_1\right).
\end{eqnarray*}
In general  
$$s_{n,N} = \left(0, 0, \ldots, 0, a_{n,1}Z_1, a_{n,1}Z_1+ a_{n+1,2}Z_2, \ldots, a_{n,1}Z_1+ a_{n+1,2}Z_2+ \ldots + a_{N,N-n+1}Z_{N-n+1}\right)$$ 
(there are $n-1$ zeros in front). 
Let $\|.\|_\infty$ denote the maximum norm in $\mathbb{H}^N$, that is, 
$$\|(y_1, y_2, \ldots, y_N)\|_\infty = \sup(\|y_1\|, \|y_2\|, \ldots, \|y_N\|).$$
Then
 \begin{eqnarray*}
\|(X_1, X_1+X_2, \ldots, X_1+X_2+\ldots+X_N)\|_\infty & = & \|s_{1,N} + s_{2,N} + \ldots+ s_{N,N}\|_\infty\\
&\leq & \|s_{1,N}\|_\infty + \ldots + \|s_{N,N}\|_\infty
 \end{eqnarray*}
That is,
   $$\sup_{1\leq n \leq N}\|X_1 + X_2 + \ldots + X_n\| \leq \sum_{n=1}^N\|s_{n,N}\|_\infty$$ and hence
       $$\mathbb{E}\left(\sup_{1\leq n \leq N}\|X_1 + X_2 + \ldots + X_n\|\right) \leq \sum_{n=1}^N \mathbb{E}(\|s_{n,N}\|_\infty).$$ 
(a) Since all the $(Z_n)$, $1\leq n \leq N$,  are independent and symmetrical, then we can apply  the classical L\'evy's inequality to each quantity $\mathbb{E}(\|s_{n,N}\|_{\infty})$ to obtain:
\begin{eqnarray*}
\mathbb{E}(\|s_{n,N}\|_{\infty}) & = & \mathbb{E}\left(\sup_{1\leq k \leq N-n} \|a_{n,1}Z_1+ a_{n+1,2}Z_2 + \ldots a_{n+k,k}Z_k\|\right), \\
 &\leq & 2 \mathbb{E}(\|a_{n,1}Z_1 + a_{n+1,2}Z_2 + \ldots+a_{N,N-n+1}Z_{N-n+1}\|)\\
                          & \leq & 2 \left(\sum_{k=1}^{N-n+1} |a_{n+k-1, k}|^2\right)^{1/2}.
\end{eqnarray*}      
Adding all these inequalities for $n = 1,2, \ldots, N$ yields
$$\mathbb{E}\left(\sup_{1\leq n \leq N}\|X_1 + \ldots + X_n\|\right) \leq 2 \sum_{n=1}^N\left(\sum_{k=1}^{N-n+1} |a_{n+k-1, k}|^2\right)^{1/2}.$$
This completes the proof of  Lemma \ref{lemma11}.\\
(b) For Lemma \ref{lemmamart}, we make use of the important observation that each sequence $s_{n,N}$ is a (discrete) martingale. This is due to the fact that each $Z_{k}$ is independent of $\sigma(Z_1, Z_2, \ldots, Z_{k-1})$ while  its coefficient $a_{n+k,k}$ is $\sigma(Z_1, Z_2, \ldots, Z_{k-1})$-measurable.  
Then in particular the norms of the elements of $s_{n,N}$ 
          $$\left(\|a_{n,1}Z_1\|, \|a_{n,1}Z_1+ a_{n+1,2}Z_2\|, \ldots, \|a_{n,1}Z_1+ a_{n+1,2}Z_2+ \ldots + a_{N,N-n+1}Z_{N-n+1}\|\right)$$
constitute a submartingale. We can therefore make use of Doob's martingale inequality: for any nonnegative submartingale $(M_k)$:
$$\mathbb{E}\left(\sup_{1\leq j \leq k} M_j^2\right) \leq 4 \mathbb{E}(M_k^2).$$  This implies
$$\mathbb{E}\left(\sup_{1\leq j \leq k} M_j\right) \leq \left(\mathbb{E}\left(\sup_{1\leq j \leq k} M_j^2\right)\right)^{1/2} \leq 2 \left(\mathbb{E}(M_k^2)\right)^{1/2}.$$
Then we obtain:
 \begin{eqnarray*}
\mathbb{E}(\|s_{n,N}\|_\infty) & = & \mathbb{E}\left(\sup_{1\leq k \leq N-n} \|a_{n,1}Z_1+ a_{n+1,2}Z_2 + \ldots a_{n+k,k}Z_k\|\right)\\
&\leq & 2 \left(\mathbb{E}(\|a_{n,1}Z_1 + a_{n+1,2}Z_2 + \ldots+a_{N,N-n+1}Z_{N-n+1}\|^2)\right)^{1/2}\\
                          & = & 2 \left(\sum_{k=1}^{N-n+1} \mathbb{E}\left(|a_{n+k-1, k}|^2\right)\right)^{1/2}.
\end{eqnarray*}
Adding these inequalities for $n=1,2,\ldots, N$ yields
$$\mathbb{E}\left(\sup_{1\leq n \leq N}\|X_1 + \ldots + X_n\|\right) \leq 2 \sum_{n=1}^N\left(\sum_{k=1}^{N-n+1} \mathbb{E}(|a_{n+k-1, k}|^2)\right)^{1/2}.$$
This concludes the proof of Lemma \ref{lemmamart}.

\subsection{Proof of Theorem \ref{th01} and Theorem \ref{thmart}} \label{sdwe3rwsa}
For any $m\geq 1$ and $N\geq 1$, we shall first estimate the quantity $$\mathbb{E}\left(\sup_{1\leq \ell \leq m}\|X_{N+1} + \ldots + X_{N+\ell}\|\right).$$
Write
\begin{eqnarray*}
X_{N+1} & = &a_{N+1,1} Z_1  + \ldots+ a_{N+1,N+1}Z_{N+1},\\
X_{N+2} &= &a_{N+2,1} Z_1   + \ldots+ a_{N+2,N+1}Z_{N+1} + a_{N+2,N+2}Z_{N+2},\\
\vdots \\
X_{N+m} &= &a_{N+m,1} Z_1  + \ldots+ a_{N+m,N+1}Z_{N+1} + a_{N+m,N+2}Z_{N+2} + \ldots+ a_{N+m,N+m}Z_{N+m}.                                                  
\end{eqnarray*}
Using the same notations as in the proof of the lemmas,  write
       $$(X_1, X_1+ X_2, \ldots, X_1+X_2 + \ldots+X_{N+m}) = \sum_{n=1}^{N+m}s_{n,N+m}.$$
Then an easy computation yields
     $$(X_{N+1}, X_{N+1}+X_{N+2}, \ldots, X_{N+1} + X_{N+2}+ \ldots+ X_{N+m}) = \sum_{n=1}^{N+m}t_{n,N+m} $$      
 where $t_{n, N+m}$ is the same as $s_{n, N+m}$ except that all the coefficients $a_{k, j}$ for $k \leq N$ and $j\leq N$ are taken to be 0. That is,
\begin{eqnarray*}
t_{n, N+m} & = & (0, 0, \ldots, 0,  b_{n,1}Z_1, b_{n,1}Z_1+ b_{n+1,2}Z_2, \ldots,\\ 
&& \qquad \qquad \qquad b_{n,1}Z_1+ b_{n+1,2}Z_2+ \ldots + b_{N+m,N+m-n+1}Z_{N+m-n+1})
\end{eqnarray*}
($t_{n, N+m}$ is a $N+m$-tuple) where 
    $$b_{k,j} = a_{k,j} \mbox{ if } k\geq N+1 \mbox{ and } b_{k,j} = 0 \mbox{ otherwise}.$$   
Then similar calculations as in the proof of the lemmas yield
\begin{eqnarray*}
\mathbb{E}\left(\sup_{1\leq \ell \leq m}\|X_{N+1} + \ldots + X_{N+\ell}\|\right) \leq \sum_{n=1}^{N+m} \mathbb{E}(\|t_{n, N+m}\|_{\infty}).
\end{eqnarray*}
(a) Let us first assume that all the coefficients $(a_{k,j})$ are constant complex numbers (non-random).   
Then for each $n \leq N+m$, L\'evy's inequality yields
\begin{eqnarray*}
\mathbb{E}(\|t_{n, N+m}\|_\infty) &\leq & 2 \mathbb{E}(\|b_{n,1}Z_1+ b_{n+1,2}Z_2+ \ldots + b_{N+m,N+m-n+1}Z_{N+m-n+1})\|)\\
& \leq & 2 \left(\sum_{k=1}^{N+m-n+1} |b_{n+k-1, k}|^2\right)^{1/2}.
\end{eqnarray*}
Hence
\begin{eqnarray*}
\mathbb{E}\left(\sup_{1\leq \ell \leq m}\|X_{N+1} + \ldots + X_{N+\ell}\|\right) \leq 2 \sum_{n=1}^{N+m} \left(\sum_{k=1}^{N+m-n+1} |b_{n+k-1, k}|^2\right)^{1/2}.
\end{eqnarray*}
Since  $$b_{k,j} = a_{k,j} \mbox{ for } k\geq N+1 \mbox{ and } b_{k,j} = 0 \mbox{ otherwise },$$ then this implies 
\begin{eqnarray*}
\mathbb{E}\left(\sup_{1\leq \ell \leq m}\|X_{N+1} + \ldots + X_{N+\ell}\|\right) &\leq& 2 \sum_{n=1}^{N+m} \left(\sum_{k= \max\{1, N+2-n\}}^{N+m-n+1} |a_{n+k-1, k}|^2\right)^{1/2}\\
& = & 2 \sum_{n=1}^{N} \left(\sum_{k= N+2-n}^{N+m-n+1} |a_{n+k-1, k}|^2\right)^{1/2}\\
&& + 2 \sum_{n=N+1}^{N+m} \left(\sum_{k= 1}^{N+m-n+1} |a_{n+k-1, k}|^2\right)^{1/2}.
\end{eqnarray*}
In particular (since all the involved terms are positive),
\begin{eqnarray}\label{df1er45sdw}
\mathbb{E}\left(\sup_{1\leq \ell \leq m}\|X_{N+1} + \ldots + X_{N+\ell}\|\right) & \leq  &
2 \sum_{n=1}^{N} \left(\sum_{k= N+2-n}^{\infty} |a_{n+k-1, k}|^2\right)^{1/2}\nonumber\\
&& + 2 \sum_{n=N+1}^{\infty} \left(\sum_{k= 1}^{\infty} |a_{n+k-1, k}|^2\right)^{1/2}.
\end{eqnarray}
Taking the limit for $m \to \infty$ gives,     
\begin{eqnarray*}
\lim_{m \to \infty} \mathbb{E}\left(\sup_{1\leq \ell \leq m}\|X_{N+1} + \ldots + X_{N+\ell}\|\right) & \leq  &
2 \sum_{n=1}^{N} \left(\sum_{k= N+2-n}^{\infty} |a_{n+k-1, k}|^2\right)^{1/2}\\
&& + 2 \sum_{n=N+1}^{\infty} \left(\sum_{k= 1}^{\infty} |a_{n+k-1, k}|^2\right)^{1/2}.
\end{eqnarray*}
Hence (by the monotone convergence theorem),
\begin{eqnarray} \label{sd23ergbg}
\mathbb{E}\left(\sup_{1\leq \ell < \infty}\|X_{N+1} + \ldots + X_{N+\ell}\|\right) & \leq  &
2 \sum_{n=1}^{N} \left(\sum_{k= N+2-n}^{\infty} |a_{n+k-1, k}|^2\right)^{1/2}\nonumber\\
&& + 2 \sum_{n=N+1}^{\infty} \left(\sum_{k= 1}^{\infty} |a_{n+k-1, k}|^2\right)^{1/2}.
\end{eqnarray}
At this stage it is important to see that the convergence of the series in (\ref{sder434e}) implies that the right hand side in (\ref{sd23ergbg}) decays to 0 as $N\to \infty$.
 Indeed,       
Write 
\begin{eqnarray*}
A_N & = & \sum_{n=1}^{N} \left(\sum_{k= N+2-n}^{\infty} |a_{n+k-1, k}|^2\right)^{1/2}, \\
B_N & = & \sum_{n=N+1}^{\infty} \left(\sum_{k= 1}^{\infty} |a_{n+k-1, k}|^2\right)^{1/2}.
\end{eqnarray*}
Clearly, 
     $$B_N = \sum_{n=N+1}^{\infty} \|d_n\|_{\ell^2}$$
(where $d_n = (a_{n,1}, a_{n+1, 2}, a_{n+2,3},\ldots)$ is the $n$-th diagonal of the matrix $(a_{n,k})$). 
Since the series $\sum_{n=1}^{\infty} \|d_n\|_{\ell^2} < \infty$, then the residual sum $\sum_{n=N+1}^{\infty} \|d_n\|_{\ell^2}$ converges to 0 as $N \to \infty$. 
Also
     $$A_N = \sum_{n=1}^N \|d_{n,N}\|_{\ell^2}$$
where $d_{n,j}$ is obtained from $d_n$ by deleting its first $j+1-n$ elements whenever $j\geq n$ and $d_{n,j}$ is just the same as $d_n$ for $j < n$, that is,
  $$d_{n,j} = (a_{j+1, j-n+2}, a_{j+2, j-n+3}, a_{j+3, j-n+3},\ldots) \mbox{ for } j\geq n.$$
Write 
   $$T_{k, N} = \sum_{n=1}^{N}\|d_{n, N+k}\|_{\ell^2},\,\,\,k=1,2,\ldots$$
Then\begin{eqnarray*}
 A_{N+k} & = & \sum_{n=1}^{N+k} \|d_{n,N+k}\|_{\ell^2} =  \sum_{n=1}^N \|d_{n,N+k}\|_{\ell^2} + \sum_{n=N+1}^{N+k} \|d_{n,N+k}\|_{\ell^2} \\
        & \leq & T_{k,N} + \sum_{n=N+1}^{N+k} \|d_{n}\|_{\ell^2} \\
        & \leq & T_{k,N} + B_N.
\end{eqnarray*}
 It is obvious that for $N$ fixed,
$$\lim_{k\to \infty} T_{k,N} = 0.$$
  Therefore,
    $$\limsup_{k \to \infty} A_{N+k} \leq B_N.$$  
Since $B_N \to 0$ for $N\to \infty$, this implies that 
$$\lim_{N\to \infty} A_N = 0.$$         
Now coming back to (\ref{sd23ergbg}) and taking the limit for $N \to \infty$ yields
\begin{eqnarray} \label{ew1sdffere}
\lim_{N\to \infty} \mathbb{E}\left(\sup_{1\leq \ell < \infty}\|X_{N+1} + \ldots + X_{N+\ell}\|\right) =0.
\end{eqnarray}
It is now an easy matter to show that this yields the almost sure convergence of the series $\sum_{n=1}^\infty X_n$. Indeed, fix $r>0$ and write for all $N \in \mathbb{N}$,
\begin{eqnarray*}
        \mathbb{P}\left(\sup_{\ell \geq 1}\|X_{N+1} + \ldots + X_{N+\ell}\|^{1/2} >\sqrt{r}\right) \leq \frac{1}{r} \mathbb{E}\left(\sup_{\ell\geq 1}\|X_{N+1} + \ldots + X_{N+\ell}\|\right).
        \end{eqnarray*} 
Then by (\ref{ew1sdffere})),
 \begin{eqnarray*}
       \lim_{N \to \infty} \mathbb{P}\left(\left(\sup_{\ell \geq 1}\|X_{N+1} + \ldots + X_{N+\ell}\|^{1/2}\right) >\sqrt{r}\right) = 0.
        \end{eqnarray*}
Then Fatou's lemma implies
     \begin{eqnarray*}
       \mathbb{P}\left(\left(\liminf_{N\to \infty} \sup_{\ell \geq 1}\|X_{N+1} + \ldots + X_{N+\ell}\|^{1/2}\right) >\sqrt{r}\right) = 0.
        \end{eqnarray*}
Since $r$ can be taken arbitrary small, this implies
     \begin{eqnarray*}
       \mathbb{P}\left(\liminf_{N\to \infty}  \sup_{\ell \geq 1}\|X_{N+1} + \ldots + X_{N+\ell}\|^{1/2} >0 \right) = 0. \end{eqnarray*}  
Therefore, almost surely,
       $$\liminf_{N\to \infty}  \sup_{\ell \geq 1}\|X_{N+1} + \ldots + X_{N+\ell}\| = 0.$$
Hence, in particular, almost surely, for each $\epsilon >0$, there exists an integer $N(\epsilon) >0$ such that 
     $$\sup_{\ell \geq 1}\|X_{N(\epsilon)+1} + \ldots + X_{N(\epsilon)+\ell}\| \leq \epsilon.$$
That is, almost surely, for each $\epsilon >0$, there exists an integer $N(\epsilon) >0$ such that for all integers $\ell\geq 1$, 
     $$\|X_{N(\epsilon)+1} + \ldots + X_{N(\epsilon)+\ell}\| \leq \epsilon$$ or equivalently
        $$\|S_{N(\epsilon)+\ell} - S_{N(\epsilon)}\| \leq \epsilon$$ (where
$S_n = X_1+ X_2+ \ldots+ X_n$  are the partial sums).       
This implies  that the sequence $(S_n)$ converges in norm almost surely. Hence the series $\sum_{n=1}^\infty X_n$ converges almost surely.  This concludes the proof of  Theorem \ref{th01}. \\
(b) In view of Lemma \ref{lemmamart}, the same argument applies for Theorem~\ref{thmart} by just replacing the quantity $|a_{k,j}|^2$ by $\mathbb{E}(|a_{k,j}|^2)$ (for all $k,j$) wherever it appears.  \hfill \qed

\subsection{Proof of Theorem \ref{eded234231s}}
We make use of the well-known Egorov theorem. Assume that the series $\sum_{n=1}^\infty X_n$ converges almost surely in $\mathbb{H}$ and consider the partial sums
 $S_n = \sum_{k=1}^n X_k$, $n \geq 1$.
Then by the Egorov theorem, for any $\delta >0$, there exists a measurable subset $A \subset \Omega$ with $\mathbb{P}(A) > 1-\delta$ such that the sequence of partial sums $(S_n(\omega))$ converges uniformly for all $\omega \in A$. That is, 
           $$\lim_{m \to \infty} \|S_{m+j}(\omega) - S_m(\omega)\| = 0 \mbox{ uniformly in } j = 1,2,\ldots$$
We may assume that $\mathbb{P}(A) > 1-K^{-1}$ where $K>0$ is the constant satisfying (\ref{cond2}). Then for all $\omega \in A$ and for  any $\epsilon>0$, there exists $m_0$ such that for all $m > m_0$ and for all $j\geq 1$, 
               $$\|S_{m+j}(\omega) - S_m(\omega)\| < \epsilon.$$
Then in particular
       \begin{eqnarray} \label{desfee34wssz}
 \mathbb{E}(\|S_{m+j} - S_m\|^2 1_{A}) \leq \epsilon^2\ \mathbb{P}(A).
 \end{eqnarray}
Write 
\begin{eqnarray*}
\mathbb{E}(\|S_{m+j}(\omega) - S_m(\omega)\|^2 1_{A}) & = & \int_A \|S_{m+j}(\omega) - S_m(\omega)\|^2  d\mathbb{P}(\omega)\\
& =  &  \int_{\Omega} \|S_{m+j}(\omega) - S_m(\omega)\|^2  d\mathbb{P}(\omega) \\
&& -  \int_{\Omega\setminus A} \|S_{m+j}(\omega) - S_m(\omega)\|^2  d\mathbb{P}(\omega) \\
& = & \mathbb{E}(\|S_{m+j} - S_m\|^2) - \mathbb{E}(\|S_{m+j} - S_m\|^2 1_{\Omega\setminus A}).  
\end{eqnarray*}
By the Cauchy-Schwarz inequality,
\begin{eqnarray*}
\mathbb{E}(\|S_{m+j} - S_m\|^2 1_{\Omega\setminus A}) & \leq & \left(\mathbb{E}(\|S_{m+j} - S_m\|^4)\right)^{1/2} 
\left(\mathbb{E}\left(1_{\Omega\setminus A}\right)\right)^{1/2}.
\end{eqnarray*}
Since condition (\ref{cond2}) implies that $$ \left(\mathbb{E}(\|S_{m+j} - S_m\|^4)\right)^{1/2} \leq  K^{1/2}\, \mathbb{E}(\|S_{m+j} - S_m\|^2),$$ it follows that
  $$\mathbb{E}(\|S_{m+j} - S_m\|^2 1_{\Omega\setminus A})  \leq  K^{1/2}\, \mathbb{E}(\|S_{m+j} - S_m\|^2) (1- P(A))^{1/2}.$$
Therefore
    \begin{eqnarray*}       
\mathbb{E}(\|S_{m+j} - S_m\|^2 1_{A}) & \geq  & \mathbb{E}(\|S_{m+j} - S_m\|^2) - K^{1/2}\, \mathbb{E}(\|S_{m+j} - S_m\|^2) (1- P(A))^{1/2}\\
& = & \left(1- K^{1/2}(1- P(A))^{1/2}\right) \mathbb{E}(\|S_{m+j} - S_m\|^2).
\end{eqnarray*}
Hence since $\mathbb{P}(A) > 1-K^{-1}$, it follows from (\ref{desfee34wssz}) that 
   $$\mathbb{E}(\|S_{m+j} - S_m\|^2) \leq \frac{\mathbb{P}(A)}{1 - \sqrt{K(1-\mathbb{P}(A))}}\, \epsilon ^2$$ uniformly for all $j\geq 1$. Then
   $$\sup_{j\geq 1}\mathbb{E}(\|S_{m+j} - S_m\|^2) \leq \frac{\mathbb{P}(A)}{1 - \sqrt{K(1-\mathbb{P}(A))}}\, \epsilon ^2.$$
  Therefore
 $$\lim_{m\to \infty} \sup_{j \geq 1} \mathbb{E}(\|S_{m+j} - S_m\|^2) = 0$$ which yields the $L^2$--convergence of  the series $\sum_{n=1}^\infty X_n$. This concludes the proof. \hfill \qed

\subsection{Proof of Theorem \ref{stopth}}
For each fixed $r>0$, consider the stopping time $\tau_r$ defined by 
$$\tau_r(\omega) = \min\{n\geq 1: \|X_1+ X_2 + \ldots + X_n\|>r\},\,\,\,\min \emptyset = \infty$$ 
(that is, the first time the random sequence $(X_1 + X_2 + \ldots+X_n, n\geq 1)$ crosses level $r$) and the associated random change of sign $(\zeta_n)$ given by  $\zeta_n(\omega) = 1$ if $\tau(\omega) \geq n$ and $\zeta_n(\omega) = -1$ if $\tau(\omega) \leq  n-1$. Then by the series  $\sum_{n=1}^\infty \zeta_n X_n$ converges in $L^2$. 
The proof of the theorem is based on the inequality:
\begin{eqnarray} \label{safawdw23}
\mathbb{P}\left(\sup_{1\leq n \leq N} \|X_1 + X_2 + \ldots+ X_n\| > r\right) & \leq  & \mathbb{P}(\|X_1 + X_2 + \ldots+ X_N\|>r) \nonumber\\
&& + \mathbb{P}(\|\zeta_1 X_1 + \zeta_2 X_2 + \ldots+ \zeta_N X_N\|>r).
\end{eqnarray}
This implies immediately that
  \begin{eqnarray} \label{s1af1awdw23}
\mathbb{P}\left(\sup_{1\leq n \leq N} \|X_1 + X_2 + \ldots+ X_n\| > r\right) & \leq  & \frac{1}{r^2}\mathbb{E}(\|X_1 + \ldots + X_N\|^2) \nonumber\\
&&  + \frac{1}{r^2}\mathbb{E}(\|\zeta_1 X_1 + \ldots + \zeta_N X_N\|^2). 
 \end{eqnarray}
 Hence obviously, \begin{eqnarray} \label{s1af1awdw22313}
\mathbb{P}\left(\sup_{1\leq n \leq N} \|X_1 + X_2 + \ldots+ X_n\| > r\right) & \leq  & \frac{1}{r^2}\left(\sup_{q\geq 1}\mathbb{E}(\|X_1 + \ldots + X_q\|^2)\right)\nonumber\\
&&  + \frac{1}{r^2}\left(\sup_{q\geq 1}\mathbb{E}(\|\zeta_1 X_1 + \ldots + \zeta_q X_q\|^2)\right). 
 \end{eqnarray}
 In the particular case where the random variables $(X_n)$ are independent, the random variable $\zeta_1 X_1 + \ldots + \zeta_N X_N$ has the same distribution as $X_1 + \ldots + X_n$ and (\ref{safawdw23}) is the classical Levy's inequality. 
To prove (\ref{safawdw23}), set $$A =  \sup_{1\leq n \leq N} \|X_1 + X_2 + \ldots+ X_n\| > r$$ and consider its  partition
(as in Kahane \cite[p 29]{Kahane_1985})
\begin{eqnarray*}
A_1&:& \|X_1\|> r\\
A_2&:& \|X_{1}\| \leq r, \|X_1 + X_{2}\| > r\\
\vdots \\
A_N &:& \|X_1\| \leq r,  \|X_1 + X_{2}\| \leq r, \ldots, \|X_1 +X_{2} + \ldots + X_{N}\| >  r.
\end{eqnarray*}
That is, $A_n = \{\omega\in \Omega: \tau_r(\omega) = n\}$. Clearly $\mathbb{P}(A) = \sum_{n=1}^N \mathbb{P}(A_n)$. Let $S_N = X_1 + X_2 + \ldots + X_N$ and $T_N = \zeta_1 X_1 + \zeta_2 X_2 + \ldots+\zeta_N X_N$.  
 For each $n$, it is clear that $\omega \in A_n$ implies that
$\|S_N(\omega)\| >  r$ or $\|T_N(\omega)\| >  r$. It is so because by definition $\omega \in A_n$ is equivalent to $\tau_r(\omega) = n$ and then
$S_N(\omega) + T_N(\omega) = 2 (X_1(\omega) + \ldots +  X_n(\omega).$  Therefore
 $$\|S_N(\omega) + T_N(\omega)\| = 2\|X_1(\omega) + \ldots +  X_n(\omega)\| > 2 r$$ from which it follows that $\|S_N(\omega)\| >  r$ or $\|T_N(\omega)\| >  r$.
Then
$$\mathbb{P}(A_n) \leq \mathbb{P}(\omega\in A_n: \|S_N(\omega)\|> r) + \mathbb{P}(\omega\in A_n: \|T_N(\omega)\|> r).$$
Adding all these inequalities for $n=1,2,\ldots, N$ yields
        $$\mathbb{P}(A) \leq \mathbb{P}(\|S_N(\omega)\|> r) + \mathbb{P}(\|T_N(\omega)\|> r)$$ which is relation (\ref{safawdw23}).
We now want to show that inequality (\ref{s1af1awdw23}) with the $L^2$--convergence of the two series $\sum_{n=1}^\infty \zeta_n X_n$ and $\sum_{n=1}^\infty X_n$ yields the almost surely convergence of the series $\sum_{n=1}^\infty X_n$. Note that the hypothesis of the theorem implies that the series $\sum_{n=1}^\infty X_n$ is itself convergent in $L^2$ as it can be seen by taking the obvious stopping time $\tau = 0$. 
Since the series $\sum_{n=1}^\infty X_n$ and $\sum_{n=1}^\infty \zeta_n X_n$ converge in $L^2$, then for any fixed $N, m \in \mathbb{N}$, (\ref{s1af1awdw22313}) yields 
 \begin{eqnarray*}
        \mathbb{P}\left(\sup_{1\leq j \leq m}\|X_{N+1} + \ldots + X_{N+j}\|>r\right) & \leq &  \frac{1}{r^2}\left(\sup_{q\geq 1}\mathbb{E}(\|X_{N+1} + \ldots + X_{N+q}\|^2)\right)\nonumber\\
&&  + \frac{1}{r^2}\left(\sup_{q\geq 1}\mathbb{E}(\|\zeta_{N+1} X_{N+1} + \ldots + \zeta_{N+q} X_{N+q}\|^2)\right).
       \end{eqnarray*}
Taking $m\to \infty$ yields (by the monotone convergence theorem)
 \begin{eqnarray*} 
        \mathbb{P}\left(\sup_{j \geq 1}\|X_{N+1} + \ldots + X_{N+j}\|>r\right) & \leq &  \frac{1}{r^2}\left(\sup_{q\geq 1}\mathbb{E}(\|X_{N+1} + \ldots + X_{N+q}\|^2)\right) +\nonumber\\
&&  \frac{1}{r^2}\left(\sup_{q\geq 1}\mathbb{E}(\|\zeta_{N+1} X_{N+1} + \ldots + \zeta_{N+q} X_{N+q}\|^2)\right).
       \end{eqnarray*}
Since the series $\sum_{n=1}^\infty X_n$ and $\sum_{n=1}^\infty \zeta_n X_n$ converge in $L^2$, then the right hand side of this inequality goes to 0 for $N \to \infty$. A similar argument as in the proof of Theorem~\ref{th01} and Theorem~\ref{thmart} yields 
that the series $\sum_{n=1}^\infty X_n$ converges almost surely.  This concludes the proof. \qed

\section{Illustrating examples}
(1)  Let $(u_n)$ be a sequence of vectors in  a Hilbert space $\mathbb{H}$ and $(X_n)$ be a sequence of real-valued random variables defined by
      $$X_n = \sum_{k=1}^n a_{n,k} Z_k$$ where  $(Z_k)$ is a sequence of independent real-valued symmetrical random variables of unit variance. By Corollary 1, if 
      \begin{eqnarray} \label{sdsde34rfs1qw}
\sum_{n=1}^\infty \left(\sum_{k=1}^\infty |a_{n+k-1,k}|^2 \|u_{n+k-1}\|^2\right)^{1/2} < \infty,
\end{eqnarray}
 then the series $\sum_{n=1}^\infty X_n u_n$ converges almost surely. 
Then assume that 
    \begin{eqnarray} \label{qssdw34erds}
a_{n,k} = O(n-k+1)^{-\alpha}\,\,\,\mbox{ for } |n-k| \to \infty,\,\,\alpha > 1 \mbox{ fixed}.
\end{eqnarray} 
Then the condition $\sum_{n=1}^\infty\|u_n\|^2 < \infty$ is sufficient for the almost sure convergence of the series $\sum_{n=1}^\infty X_n u_n$. In other words if the coefficients $a_{n,k}$ decay to $0$ as $|n-k|\to \infty$ quickly enough, then the random variables $(X_n)$ can be treated as independent random variables with respect to the convergence of the series $\sum_{n=1}^\infty X_n u_n$. Indeed, 
$$\sum_{k=1}^\infty |a_{n+k-1,k}|^2 \|u_{n+k-1}\|^2 \leq K \sum_{k=1}^\infty n^{-2\alpha} \|u_{n+k-1}\|^2 \leq K_1 n^{-2\alpha}$$ for some constants $K, K_1$ (independent of $n$). This yields (\ref{sdsde34rfs1qw}) because $\alpha > 1$. This is an important generalisation of a result which was only known for independent random variables and few particular dependent variables (negatively associated random variables as discussed in the introduction). 
 \\
If instead we assume that
      \begin{eqnarray*} 
a_{n,k} = O(n-k+1)^{-\alpha}\,\,\,\mbox{ for } |n-k| \to \infty,\,\,\alpha \leq  1 \mbox{ fixed},
\end{eqnarray*}
the convergence of the numerical series $\sum_{n=1}^\infty \|u_n\|^2 < \infty$ may not be enough for the almost sure convergence of $\sum_{n=1}^\infty X_n$. However a faster decay such as
  $$\|u_n\|^2 = O\left(n^{-(1+ 2 \beta)}\right)\,\,\mbox{ for } n\to \infty \mbox{ with } \beta > 1-\alpha  $$
guarantees the almost sure convergence of the series $\sum_{n=1}^\infty X_n u_n$ and its (non-random) change of signs series $\sum_{n=1}^\infty \pm X_n u_n$.  \\
(2) Let us consider in particular the standard trigonometric series $\sum_{n=1}^\infty f_n(x)$, where
 $$f_n(x) = A_n \cos n x + B_n \sin nx,\,\,A_n, B_n\in \mathbb{R},\, x\in [0, 2\pi]$$ and the corresponding random series
$\sum_{n=0}^\infty X_n f_n(x)$  where $(X_n)$ is a sequence of symmetrical  real-valued random variables with finite variance.
 In  $L^2[0, 2\pi]$,  $\|f_n\|^2 = A_n^2 + B_n^2.$              
 In the case where the random variables $(X_n)$ are independent (and have unit variance), it is well-known that the series  $\sum_{n=0}^\infty X_n f_n(x)$ converges almost surely in $L^2[0, 2\pi]$  if and only if $\sum_{n=1}^\infty (A_n^2 + B_n^2) < \infty.$ (See for example the book by Zygmund \cite[p 214]{Zygmund}.)            
The discussion above shows that the condition  $\sum_{n=1}^\infty (A_n^2 + B_n^2) < \infty$  is also sufficient for the almost sure convergence  of the series $\sum_{n=0}^\infty X_n f_n(x)$ in $L^2[0, 2\pi]$ when each $X_n$ can be written as a linear combination of independent symmetrical random variables of unit variances 
$X_n = \sum_{n=1}^n a_{n,k} Z_n$ such that (\ref{qssdw34erds}) holds. \\
(3) Consider the classical fractional Gaussian noise of index $0 \leq H < 1$. It is a sequence of Gaussian real variables $(\Delta_n)$ defined on the same probability space with the covariance structure:
     \begin{eqnarray} \label{ew32wswaws2}
\mathbb{E}(\Delta_n \Delta_m) = {\scriptstyle\frac{1}{2}}|n-m+1|^{2H} + {\scriptstyle\frac{1}{2}}|n-m-1|^{2H} -|n-m|^{2H}. 
\end{eqnarray}  
Alternatively for $H\ne 0$,  $\Delta_n = B(n+1)-B(n)$ where $\{B(t): t\geq 0\}$ is the classical fractional Brownian motion process of Hurst index $H$. (See for example the book by Nourdin \cite{Nourdin_Ivan} for some details on this process.)
A model of such sequence can be obtained by considering a sequence $(Z_n)$ is independent and identically distributed standard real Gaussian variables and write
  $$\Delta_n = \alpha_{n,1} Z_1 +\ldots + \alpha_{n,n}Z_n $$ where $\alpha_{n,k}$ are real numbers satisfying the following relations (obtained from the covariance structure (\ref{ew32wswaws2})):
  \begin{eqnarray*}
  \sum_{j=1}^k \alpha_{n,j} \alpha_{k,j} & = & \mathbb{E}(\Delta_n \Delta_k),\,\,\,\,k \leq n \\
\sum_{k=1}^n \alpha_{n,j}^2 & = & 1.
 \end{eqnarray*}
Functionals of the form 
   $$\sum_{n=1}^\infty f(n) (B(n+1) - B(n)) = \sum_{n=1}^\infty f(n) \Delta_n$$ for a function $f:\mathbb{R} \to \mathbb{R}$ which is non-random (or random but adapted to the process $\{B(t): t\geq 0\}$) are closely related to  stochastic integration.  It is interesting to know under which conditions on the function $f$, is such series almost surely convergent.  For the classical Brownian motion (corresponding to $H = \frac{1}{2}$), the increments $(\Delta_n)$ are independent and hence, as already discussed, the condition $\sum_{n=1}^\infty \mathbb{E}(\|f(n)\|^2)$ is sufficient for the convergence of the series $\sum_{n=1}^\infty f(n) \Delta_n$. In the general case, in view of Corollary 1, if
        \begin{eqnarray*} \label{sdsde34rfs}
\sum_{n=1}^\infty \left(\sum_{k=1}^\infty |\alpha_{n+k-1,k}|^2\ \mathbb{E}(\|f{n+k-1}\|^2)\right)^{1/2} < \infty,
\end{eqnarray*}
then the series $\sum_{n=1}^\infty f(n) \Delta_n$ also converges almost surely. This can be generalised to other Gaussian processes.

 \paragraph{Remark:} An anonymous referee asked under which conditions the re-arrangement of the terms of the series $\sum_{n=1}^\infty X_n$ is possible so that equality (\ref{s423de234sdw}) holds, that is:
  \begin{eqnarray} \label{dsdedsd34f11} 
 \sum_{n=1}^\infty X_n = \sum_{k=1}^\infty \left(\sum_{n=1}^\infty a_{n, k}\right) Z_k.
\end{eqnarray} 
 A general condition for the almost sure convergence of all possible re-arrangements of a random series $\sum_{n=1}^\infty \xi_n$ in a Banach space $B$ (such series are said to converge almost surely unconditionally) is given by the following result from Kvaratskhelia \cite{Kvaratskhelia} (see Proposition 2.2.3 page 245 and Theorem 2.2.5 page 248).
 \begin{proposition}
Let $\sum_{k=1}^\infty \xi_k$ be an arbitrary random series in a Banach space  $B$. Then the series
$\sum_{k=1}^\infty \xi_k$  converges almost surely unconditionally in  $B$ if 
  $$\mathbb{E}\left(\sup_{\|x^*\|_{B^*} \leq 1} \sum_{k=1}^\infty \left|\langle x^*, \xi_k\rangle\right|\right) < \infty$$ and 
  $$\lim_{n \to \infty} \mathbb{E}\left(\sup_{\|x^*\|_{B^*} \leq 1} \sum_{k=n}^\infty \left|\langle x^*, \xi_k\rangle\right|\right) = 0$$
  where $B^*$ is the dual of $B$ and for $x^* \in B^*$, 
         $$\|x^*\| = \sup\{|\langle x^*, x\rangle|: x \in B \mbox{ and } \|x\| \leq 1\}.$$ 
  For Gaussian random variables $(\xi_k)$ these conditions are also necessary for the series $\sum_{k=1}^\infty \xi_k$ to converge almost surely unconditionally. 
\end{proposition}
 In the particular case of numerical random series $(B = \mathbb{R})$, these conditions reduce to 
               $$\mathbb{E}\left(\sum_{k=1}^\infty \left|\xi_k\right|\right) < \infty.$$
In our case, the unconditional almost sure convergence of the series
 $$\sum_{n=1}^{\infty} \sum_{k=1}^n a_{n,k} Z_k$$ is guaranteed (in the case particular case the underlying Hilbert space $\mathbb{H}$ is taken as $\mathbb{R}$) 
  by the conditions
              $$\mathbb{E}\left(\sum_{n=1}^\infty \sum_{k=1}^n \left|a_{n,k}| |Z_k\right|\right) < \infty.$$
 For standard Gaussian random variables $(Z_n)$, this is equivalent to 
             \begin{eqnarray} \label{dsds3ereds} \sum_{n=1}^\infty \sum_{k=1}^n \left|a_{n,k}\right| < \infty. 
             \end{eqnarray}
As it can be readily seen, this condition is much stronger than our condition 
                \begin{eqnarray*}
\sum_{n=1}^\infty \left(\sum_{k=1}^\infty |a_{n+k-1,k}|^2\right)^{1/2} < \infty,
\end{eqnarray*}
for the almost sure convergence of the initial series 
         $$\sum_{n=1}^\infty \left(\sum_{k=1}^n a_{n,k} Z_k\right).$$ 
One can immediately see that  condition (\ref{dsds3ereds}) is not necessary for equality  (\ref{dsdedsd34f11}) to hold by simply taking $a_{n,k} = 0$ for $k \ne n$, (that is, the random variables $X_n = \sum_{k=1}^n a_{n,k} Z_k$, $n \geq 1$, are independent) and $$\sum_{n=1}^\infty |a_{n,n}|^2 < \infty \mbox{  but } \sum_{n=1}^\infty |a_{n,n}| = \infty.$$ 
 The reviewer's question can be extended to the following: 
   Under which necessary and sufficient condition does the equality
          \begin{eqnarray*} 
\sum_{n=1}^\infty \left(\sum_{k=1}^n a_{n,k} Z_k\right) = \sum_{k=1}^\infty \left(\sum_{n=1}^\infty a_{n, k}\right) Z_k
\end{eqnarray*} 
hold? This question requires further investigations.

\section{Concluding remarks}
The paper deals with random series $\sum_{n=1}^\infty X_n$ where each term $X_n$ is a linear combination $\sum_{k=1}^n a_{n,k} Z_k$ where $(a_{n,k})$ are constant complex numbers and $(Z_n)$ is a sequence of independent symmetrical random variables of unit variance in a Hilbert space $\mathbb{H}$. We have obtained an  explicit  sufficient condition on the coefficient matrix $(a_{n,k})$ for the almost sure convergence of the series $\sum_{n=1}^\infty X_n$ and all its (non-random) change of signs $\sum_{n=1}^\infty \pm X_n$, namely:
\begin{eqnarray*} 
\sum_{n=1}^\infty \left(\sum_{k=1}^\infty \|a_{n+k-1,k}\|^2\right)^{1/2} < \infty.
\end{eqnarray*}
 In the proof we have used an extension of the classical L\'evy  inequality which in itself has some independent interest. This result is extended to the case where the coefficients $a_{n,k}$ are themselves assumed to be random variables taking values in the complex plane, provided that they satisfy a  natural measurability condition: each $a_{n,k}$ is measurable with respect to the $\sigma$-algebra spanned by the past $Z_1, Z_2, \ldots, Z_{k-1}$. In that case, the condition becomes 
  \begin{eqnarray*} 
\sum_{n=1}^\infty \left(\sum_{k=1}^\infty \mathbb{E}(\|a_{n+k-1,k}\|^2)\right)^{1/2} < \infty.
\end{eqnarray*}
Our results cover the particular case where the random variables $X_n$ are independent (corresponding to $a_{n,k} = 0$ for $n\ne k$) and the obvious dependence case where all the $X_n$ are collinear random variables (corresponding to $a_{n,k} = 0$ for all $k\geq 2$). For the particular case of independent random variables, this condition  is also necessary for the almost sure convergence of all the series $\sum_{n=1}^\infty \pm X_n$ (under the assumption that  $\mathbb{E}(\|X_n\|^4) \leq C \mathbb{E}(\|X_n\|^2 < \infty$ for some fixed constant $C$).  It would be interesting to know whether (under the same assumption) our condition is also necessary for the almost sure convergence in the general case of the series $\sum_{n=1}^\infty \pm X_n$ considered in this paper. This necessitates further investigation. 

The paper also gives a condition under which $L^2$-convergence implies almost sure convergence for the general series $\sum_{n=1}^\infty X_n$ and finally shows that the almost sure convergence can be deduced from the $L^2$-convergence of some modified series defined by stopping times. 
We hope that this paper will be of some value for those who would like to take further the study of random series of dependent random variables.  
\paragraph{Acknowledgements} I thank the anonymous referees whose comments helped to improve the paper. 
This research received with gratitude funding  from the College of Economics and Management Sciences of the University of South Africa. Some part of this work was completed during a research visit to the African Center of Excellence in Data Science (ACE--DS) of the University of Rwanda and I thank the Center management for its hospitality.


\begin{thebibliography}{100}

\bibitem{Esary} Esary, J., Proschan, F. and Walkup, D. (1967). Association of random variables with applications. {\it Ann. Math. Stat.} 38, 1466--1474. 

\bibitem{Huan} Huan, N.V. (2015). On the complete  convergence for sequences of random vectors in Hilbert spaces. {\it Acta Math. Hungar.} 147(1), 205--219.

\bibitem{Joag-Dev} Joag-Dev, K. and Proschan, F. (1983). Negative association of random variables with applications. {\it Ann. Stat.} 11, 286--295.

\bibitem{Kahane_1985}
Kahane, J.-P. 1985.
\newblock {\em Some random series of functions, 2nd ed.}
\newblock Cambrigde University Press, Cambridge. 

\bibitem{Kvaratskhelia} Kvaratskhelia, V.V. (2014).  Unconditional convergence of functional series in problems of probability theory, {\it Journal of Mathematical Sciences}, 200(2), 143--294.

\bibitem{Ko} Ko, M.H, Kim, T.S. and Han, K.H. (2009). A note on the almost sure convergence for dependent random variables in Hilbert space. {\it Journal of Theoretical Probability}, 22, 506--513. 

\bibitem{Levental} Levental, S.,  Mandrekar, V. and Chobonyan, S.A. (2011). Towards Nikishin's theorem on the almost sure convergence of rearrangements of functional series. {\it Functional  Analysis and Applications} 45(1), 33-45.

\bibitem{Matula} {Matu\l{}a}, P. (1992).  A note on the almost sure convergence of sums of negatively dependent random variables. {\it Statisitics and Probability  Letters} 15, 209--2013.

\bibitem{Nourdin_Ivan} Nourdin, I. (2012). {\it Selected aspects of fractional Brownian motion.}  Bocconi University Press, Springer-Verlag.


\bibitem{Shiryaev}  Shiryaev, A.N. (2016). {\it Probability--1, 3rd edition}. Springer. 

\bibitem{Vakhania} Vakhania, N.N., Tarieladze, V.I. and Chobanyan, S.A. (1987). {\it Probability Distributions
on Banach Spaces}, Mathematics and its Applications, vol. 14, D. Reidel Publishing Co., Dordrecht.

\bibitem{Wu_Jiang} Wu, Q. and Jiang, Y. (2018). Some limiting behavior for asymptotically negative associated random variables. {\it Probability in Engineering and Information Sciences}, 32, 58-66.


\bibitem{Zygmund} Zygmund, A. 1959. {\it  Trigonometric Series, Vol. I.} Cambridge University Press.
\end{thebibliography}
\end{document}